\documentclass{article}
\usepackage[utf8]{inputenc}
\usepackage{amsmath}
\usepackage{graphicx}

\title{Renormalized Reduced Order Models with Memory for Long Time Prediction}
\author{Jacob Price, Panos Stinis \\ 
University of Washington, Pacific Northwest National Laboratory}

\begin{document}
\maketitle
\begin{abstract}
We examine the challenging problem of constructing reduced models for the long time prediction of systems where there is no timescale separation between the resolved and unresolved variables. In previous work we focused on the case where there was only transfer of activity (e.g. energy, mass) from the resolved to the unresolved variables. Here we investigate the much more difficult case where there is {\it two-way} transfer of activity between the resolved and unresolved variables. Like in the case of activity drain out of the resolved variables, even if one starts with an exact formalism, like the Mori-Zwanzig (MZ) formalism, the constructed reduced models can become unstable. We show how to remedy this situation by using dynamic information from the full system to {\it renormalize} the MZ reduced models. In addition to being stabilized, the renormalized models can be accurate for very long times. We use the Korteweg-de Vries equation to illustrate the approach. The coefficients of the renormalized models exhibit rich structure, including algebraic time dependence and incomplete similarity. 
\end{abstract}

\section{Introduction}

There exist many systems whose sizes preclude the complete simulation of their dynamics given finite computational power and time. Applications such as physical chemistry, nuclear engineering, plasma physics, and climate modeling regularly produce prohibitively large systems. Reduced order modeling seeks to reduce a prohibitively large system to a computationally realizable problem size (resolution) while maintaining the essential features of the dynamics. Moreover, even if one could fully resolve a very complex system, reduced order modeling is an attractive concept because it allows one to follow the most important features of the dynamics while accounting for the rest implicitly. In both cases, the process of model reduction itself can be particularly helpful. It allows us to probe the dynamic interaction of the system constituents and often extract remarkable structure that is not at all obvious for complex systems.  

The principal goal of reduced order modeling is to reduce the system in such a way that the dynamics for the reduced set of variables agree as closely as possible with what their dynamics would be in a fully resolved system. More often than not, the degrees of freedom that are not simulated (unresolved) affect strongly the degrees of freedom that are kept. The Mori-Zwanzig (MZ) formalism \cite{zwanzig1961memory,mori1965} is an exact formalism for the reduction of the dynamics of a full system to the dynamics of a reduced set of variables. The impact of the unresolved variables on the reduced set of variables manifests in several ways, including a memory term \cite{chorin2000optimal,chorin2002optimal}. The memory term depends on the history of the trajectory of a variable up until the current time. It is fundamental to reduced order models, and usually represents a significant computational challenge. This is particularly true when there is {\it no} timescale separation between the resolved and unresolved variables. 

Due to the significant cost of computing the memory term, several approximations and expansions have been suggested \cite{chorin2002optimal,hald2007optimal,stinis2007higher,lei2016,parish2017}. In certain cases such reduced models have resulted in significant improvement of our predictive ability. However, some of them can also suffer from instabilities. Ideas from renormalization theory in physics \cite{goldenfeld1992,georgi1993} have been recently used to stabilize reduced models \cite{stinis2013renormalized,stinis2015renormalized}. The form of the problematic terms is retained but their strength is controlled by ``renormalized" coefficients. These coefficients are chosen so that the predictions of the reduced model match observed quantities from the full system. We note that dealing with time-dependent systems of equations makes the application of the renormalization concept rather delicate. 

Such renormalized reduced models have been constructed so far for systems where there is only drain of activity (energy) from the resolved to the unresolved variables \cite{stinis2013renormalized,stinis2015renormalized}. Here we extend this construction to the much more difficult case of systems with {\it two-way} transfer of activity between resolved and unresolved variables. In order to do that, first we derive a new way of approximating the memory term which we call the ``complete memory approximation''. We provide a way of automating the calculation of the terms of the approximation through the use of software e.g. Mathematica. Then we show, using the example of the Korteweg-de Vries (KdV) equation with small dispersion, that renormalization is necessary even if the amount of activity transferred back and forth between resolved and unresolved variables is small. We describe a way to renormalize the reduced models in such a case. Our construction allows us to reveal the rich structure of the renormalized coeffcients. This includes algebraic time dependence and incomplete similarity in the magnitude of the dispersion and the resolution of the reduced model. Finally, the renormalized reduced models are used successfully for the long time prediction of solutions to the KdV equation.   

The paper is organized as follows. Section 2 gives a brief overview of the Mori-Zwanzig formalism, the complete memory approximation, and the main idea behind the renormalized models. In Section 3, we construct reduced order models for the KdV equation utilizing the complete memory approximation. We also compute the renormalized coefficients for these reduced models. Section 4 contains numerical results of the long time behavior of solutions of the renormalized reduced models. The final section involves a discussion of the results and suggestions for future work.

\section{The Mori-Zwanzig Formalism}

Consider a system of (in general nonlinear) autonomous ordinary differential equations
\begin{equation}
\frac{d \mathbf{u}(t)}{dt} = \mathbf{R}(\mathbf{u})\label{orig}
\end{equation}
with an initial condition $\mathbf{u}(0)=\mathbf{u}^0$. For example, if one considers a partial differential equation, spectral or finite volume methods can be employed to convert the infinite-dimensional PDE to a system of ODEs. Let $\mathbf{u}(t) = \{u_k(t)\}$, $k\in F\cup G$. We separate $\mathbf{u}(t)$ into resolved variables $\tilde{\mathbf{u}}=\{u_k(t)\}$, $k\in F$ and unresolved variables $\hat{\mathbf{u}}=\{u_k(t)\}$, $k\in G$. Let $R_k(\mathbf{u})$ be the $k$th entry in the vector-valued function $\mathbf{R}(\mathbf{u})$. We can transform this nonlinear system of ODEs \eqref{orig} into a linear system of PDEs by way of the Liouvillian operator, also known as the generator of the Koopman operator \cite{koopman1931hamiltonian}:

\begin{equation}
\mathcal{L}=\sum_{k\in F\cup G} R_k(\mathbf{u}^0)\frac{\partial}{\partial u_k^0}.\label{L}
\end{equation}
It can be shown \cite{chorin2000optimal} that if $\phi(\mathbf{u}^0,t)$ satisfies
\begin{equation}
\frac{\partial \phi(\mathbf{u}^0,t)}{\partial t} = \mathcal{L}\phi,\;\; \phi(\mathbf{u}^0,0)=f(\mathbf{u}^0),\label{liouville}
\end{equation}
for a given function $f(\mathbf{u})$ and an initial condition $\mathbf{u}^0$, then $f(\mathbf{u}(t)) = \phi(\mathbf{u}^0,t)$. Thus, the nonlinear ODE and linear PDE are equivalent. In particular, if we consider $\phi_k$ the solution of \eqref{liouville} with $f(\mathbf{u}^0)=u_k^0$, $\phi_k(\mathbf{u}^0,t) = u_k(t)$, so the evolution of any component $u_k$ can be expressed in terms of a linear PDE. 


Consider the space of functions of $\mathbf{u}$. Let $P$ be an orthogonal projection onto the subspace of functions depending only on the resolved variables $\tilde{\mathbf{u}}$. For example, $Pf$ might be the conditional expectation of $f$ given the resolved variables. Let $Q=I-P$. If we write $\phi_k(\mathbf{u}^0,t) = e^{t\mathcal{L}}u_k^0$ and Dyson's formula:
\begin{equation}
e^{t\mathcal{L}}=e^{tQ\mathcal{L}}+\int_0^t e^{(t-s)\mathcal{L}}P\mathcal{L}e^{sQ\mathcal{L}}\,\mathrm{d}s, \label{dyson}
\end{equation}
we find
\begin{equation}
\frac{\partial}{\partial t}e^{t\mathcal{L}}u_k^0 = e^{t\mathcal{L}}P\mathcal{L}u_k^0 +e^{tQ\mathcal{L}}Q\mathcal{L}u_k^0 + \int_0^t e^{(t-s)\mathcal{L}}P\mathcal{L}e^{sQ\mathcal{L}}Q\mathcal{L}u_k^0\,\mathrm{d}s.\label{MZ}
\end{equation}This is the Mori-Zwanzig identity. It represents an alternative, but exact way of writing the full dynamics. The first term on the right hand side in \eqref{MZ} is called the Markovian term, because it depends only on the instantaneous values of the resolved variables. The second term is called `noise' and the third is called `memory'. If we project again
\begin{equation}
\frac{\partial}{\partial t}Pe^{t\mathcal{L}}u_k^0  = Pe^{t\mathcal{L}}P\mathcal{L}u_k^0+P\int_0^te^{(t-s)\mathcal{L}}P\mathcal{L}e^{sQ\mathcal{L}}Q\mathcal{L}u_k^0\,\mathrm{d}s.\label{reduced_MZ}
\end{equation}Here, we made use of the fact that
\[Pe^{tQ\mathcal{L}}Q\mathcal{L}u_k^0= P\left[I+tQ\mathcal{L}+t^2(Q\mathcal{L})^2+\dots\right]Q\mathcal{L}u_k^0 = 0
\]because $PQ=0$. For $k\in F$, \eqref{reduced_MZ} describes the projected dynamics of the resolved variables. The system is not closed, however, due to the presence of the orthogonal dynamics operator $e^{sQ\mathcal{L}}$ in the memory term. In order to simulate the dynamics of \eqref{reduced_MZ} exactly, one needs to evaluate the second term which requires the dynamics of the unresolved variables. One key fact must be understood: reducing a large system to one of comparatively fewer variables necessarily introduces a memory term encoding the interplay between the unresolved and resolved variables. Dropping both the noise and memory terms and simulating only the ``average'' dynamics (the Markovian term) may not accurately reflect the dynamics of the resolved variables in the full simulation. Any reduced dynamical model must in some way approximate or compute the memory term, or argue convincingly why the memory term is negligible. In fact, we will demonstrate that even in cases where the magnitude of the memory term is very small, it must be included in order to produce accurate simulations of the resolved modes.

\subsection{The Complete Memory Approximation}

Constructing a reduced order model (ROM) from the Mori-Zwanzig formalism requires approximating the memory integral in terms of only the resolved variables. In \cite{stinis2013renormalized}, a class of ROMs is derived under the assumption of the {\it almost} commutativity of the $P\mathcal{L}$ and $Q\mathcal{L}$ operators (see also Appendix \ref{bch_expansion}). Here, we derive a new class of ROMs that avoids this assumption. If $P\mathcal{L}Q\mathcal{L}=Q\mathcal{L}P\mathcal{L}$, the two classes are equivalent. For this reason, we call our new class the ``complete memory approximation.''

We rewrite the memory term using the definitions of $e^{-s\mathcal{L}}$ and $e^{sQ\mathcal{L}}$ and then computing the integral termwise:
\begin{align}
P\int_0^t&e^{(t-s)\mathcal{L}}P\mathcal{L}e^{sQ\mathcal{L}}Q\mathcal{L}u_k^0\,\mathrm{d}s\notag\\
=&Pe^{t\mathcal{L}}\int_0^t \left(\sum_{i=0}^{\infty}\frac{(-1)^is^i}{i!}\mathcal{L}^i\right)P\mathcal{L}\left(\sum_{j=0}^\infty \frac{s^j}{j!}(Q\mathcal{L})^j\right)Q\mathcal{L}u_k^0\,\mathrm{d}s\notag\\
=&Pe^{t\mathcal{L}}\left(\sum_{i=0}^\infty \sum_{j=0}^\infty \frac{(-1)^it^{i+j+1}}{i!j!(i+j+1)}\mathcal{L}^iP\mathcal{L}(Q\mathcal{L})^{j}Q\mathcal{L}u_k^0\right)\label{eq:novel}
\end{align}
Consider writing the first few terms arranged by powers of $t$:
\begin{align}
P\int_0^t&e^{(t-s)\mathcal{L}}P\mathcal{L}e^{sQ\mathcal{L}}Q\mathcal{L}u_k^0\,\mathrm{d}s \notag\\ =&tPe^{t\mathcal{L}}\left[P\mathcal{L}Q\mathcal{L}\right]u_k^0 + \frac{t^2}{2}Pe^{t\mathcal{L}}\left[P\mathcal{L}Q\mathcal{L}Q\mathcal{L}-\mathcal{L}P\mathcal{L}Q\mathcal{L}\right]u_k^0+O(t^3) \label{unclosed}
\end{align}The $O(t)$ term is the $t$-model which has been used to great success in \cite{bernstein2007optimal,chorin2002optimal,chorin2000optimal}. The $O(t^2)$ term, presents a new problem. The first term in it can be computed in a manner similar to the $t$-model, but the second term is not projected prior to its evolution. This makes it impossible to compute as part of a reduced-order model, as it would depend on the specific dynamics of the unresolved quantities.

To close the model in the resolved variables we start with constructing an additional reduced order model for the unclosed term. This term is $Pe^{t\mathcal{L}}\mathcal{L}P\mathcal{L}Q\mathcal{L}u_k^0$. First note that

\begin{equation}
Pe^{t\mathcal{L}}\mathcal{L}P\mathcal{L}Q\mathcal{L}u_k^0=\frac{\partial}{\partial t} Pe^{t\mathcal{L}}P\mathcal{L}Q\mathcal{L}u_k^0.
\end{equation}That is, it is the derivative of the $t$-model term itself. Now consider a reduced order model for this derivative under the Mori-Zwanzig formalism again:
\begin{align}
\frac{\partial}{\partial t} &Pe^{t\mathcal{L}}P\mathcal{L}Q\mathcal{L}u_k^0 =Pe^{t\mathcal{L}}P\mathcal{L}P\mathcal{L}Q\mathcal{L}u_k^0+P\int_0^t e^{(t-s)\mathcal{L}}P\mathcal{L}e^{sQ\mathcal{L}}Q\mathcal{L}P\mathcal{L}Q\mathcal{L}u_k^0\,\mathrm{d}s\notag\\
=&Pe^{t\mathcal{L}}P\mathcal{L}P\mathcal{L}Q\mathcal{L}u_k^0 + Pe^{t\mathcal{L}}\left(\sum_{i=0}^\infty\sum_{j=0}^\infty \frac{(-1)^it^{i+j+1}}{i!j!(i+j+1)}\mathcal{L}^iP\mathcal{L}(Q\mathcal{L})^{j+1}P\mathcal{L}Q\mathcal{L}u_k^0\right).\label{tmodeloftmodel}
\end{align}

If we replace $Pe^{t\mathcal{L}}\mathcal{L}P\mathcal{L}Q\mathcal{L}u_k^0$ in \eqref{unclosed} with \eqref{tmodeloftmodel}, only the first term of  \eqref{tmodeloftmodel} contributes at the $O(t^2)$ level. All the double sum terms contribute at $O(t^3)$ and higher. We find 
\begin{align}
P\int_0^t&e^{(t-s)\mathcal{L}}P\mathcal{L}e^{sQ\mathcal{L}}Q\mathcal{L}u_k^0\,\mathrm{d}s \notag\\ =&tPe^{t\mathcal{L}}\left[P\mathcal{L}Q\mathcal{L}\right]u_k^0 - \frac{t^2}{2}Pe^{t\mathcal{L}}\left[P\mathcal{L}P\mathcal{L}Q\mathcal{L}-P\mathcal{L}Q\mathcal{L}Q\mathcal{L}\right]u_k^0+O(t^3) \label{closed}
\end{align}
where all the terms are now projected prior to evolution and so it involves {\it only} resolved variables.

We can naturally extend this to higher orders in a straightforward manner. Any time we are left with a term that is not projected prior to evolution, we can construct a reduced-order model for that term using the MZ formalism. The result is a telescoping construction of reduced order models. It can become tedious to derive higher ordered terms by hand, but we developed a symbolic algorithm in Mathematica to compute them automatically (the software is available upon request).

\subsection{Renormalized Memory Approximations}

The complete memory approximation framework provides a series representation of the memory integral. Different ROMs can be created by truncating the series at different terms. In this case, the differential equation for a resolved mode is:
\begin{equation}
\frac{du_k}{dt} = R_k^0(\hat{\mathbf{u}}) + \sum_{i=1}^N \frac{(-1)^{i+1}t^i}{i!} R_k^i(\hat{\mathbf{u}}),\label{nonrenormalized}
\end{equation}where $R_k^0(\hat{\mathbf{u}})$ is the Markov term, $R_k^1(\hat{\mathbf{u}}) = Pe^{t\mathcal{L}}[P\mathcal{L}Q\mathcal{L}]u_k^0$ is the $t$-model, and higher ordered terms are found by grouping similar powers of $t$ in the complete memory approximation. For the examples we have tried, the resulting ROMs are unstable.

Similarly structured ROMs have been stabilized through \emph{renormalization}, in which we attach additional coefficients to each term in the series, such that the terms represent an \emph{effective memory}, given knowledge only of the resolved modes \cite{stinis2013renormalized,stinis2015renormalized}. The evolution equation for a reduced variable then becomes
\begin{equation}
\frac{du_k}{dt} = R_k^0(\hat{\mathbf{u}}) + \sum_{i=1}^N \alpha_i(t)t^i R_k^i(\hat{\mathbf{u}}).\label{renormalized}
\end{equation}Here, we allow the renormalization coefficients $\alpha_i(t)$ to be time dependent. These coefficients must be chosen in a way that captures information we know about the memory term. Designating them is done in a problem-specific manner, frequently by comparing the reduced system to a larger system prior to the point where the larger system becomes unresolved. Concrete examples and details will be provided in the following sections.

\section{Reduced Order Models of the Korteweg-de Vries Equation with Small Dispersion}\label{kdv}

The Korteweg-de Vries equation with small dispersion $\epsilon$ is
\begin{equation}
u_t + uu_x + \epsilon^2 u_{xxx} = 0. \label{kdv_equation}
\end{equation}We will consider solving this equation on $[0,2\pi]$ with periodic boundary conditions and initial condition $u^0 = u(x,0)$. 

Renormalized MZ models have already been applied to the Burgers equation which corresponds to the case $\epsilon=0$ in \eqref{kdv_equation} \cite{stinis2013renormalized,stinis2015renormalized}. The Burgers equation develops singularities in the form of shocks in finite time. For $\epsilon \neq 0,$ the dispersive term precludes a finite time singularity. Instead, the solution can be fully resolved with $O(1/\epsilon)$ Fourier modes \cite{venakides1987zero}. Additionally, the presence of a dispersive term causes energy to flow from the resolved modes to the unresolved modes \emph{and back}, unlike Burgers where there is a constant drain of energy out of the resolved variables. 

We have to make a comment here about terminology. In the dispersive equation community (which includes KdV), the square magnitude of a Fourier mode is called the ``mass'' of the mode while in the fluid dynamics community (which includes the study of Burgers) it is called the ``energy'' of the mode. Since we will consider the KdV equation we will from now on use the word ``mass.'' 

Because of the periodic boundary conditions, we use Fourier series as a basis for the solution. That is, let \[u(x,t) = \sum_{k\in F\cup G} u_k(t) e^{ikx}\]where $F\cup G = \left[ -M,\dots, M-1\right]$ and $F=\left[-N,\dots, N-1\right]$ for $N<M$. We call $F$ the \emph{resolved} modes and $G$ the \emph{unresolved} modes. Let $\mathbf{u} =\{u_k(t)\}_{k\in F\cup G}$. We partition $\mathbf{u}=(\hat{\mathbf{u}},\tilde{\mathbf{u}})$ where $\hat{\mathbf{u}} = \{u_k\}_{k\in F}$ and $\tilde{\mathbf{u}}=\{u_k\}_{k\in G}$. Our goal is to construct a reduced order model for each component $u_k(t)$ of $\hat{\mathbf{u}}$. $\mathbf{u}$ is the ``full'' model, which will be fully resolved if $G$ exceeds the maximal resolution for the chosen $\epsilon$.

The equation of motion for the Fourier mode $u_k$ is
\begin{equation}
\frac{d u_k}{dt} =R_k(\mathbf{u}) = i\epsilon^2 k^3 u_k-\frac{ik}{2}\sum_{\substack{p+q=k\\p,q\in F\cup G}}u_p u_q.\label{full}
\end{equation}
The convolution sum in the second term on the right hand side can be computed efficiently by transforming data to real space and computing the sum there as the FFT of the product of the real space solution.

We define $\mathcal{L}$ as in \eqref{L} such that $\mathcal{L}u_k^0 = R_k(\mathbf{u}^0)$. We must define the projection operator $P$. Consider a function $h(\mathbf{u})$ that depends on all the Fourier modes. We define $Ph(\mathbf{u}^0)=Ph(\hat{\mathbf{u}}^0,\tilde{\mathbf{u}}^0)=h(\hat{\mathbf{u}}^0,0)$. That is, we set each unresolved variable to zero. In order to remain consistent with our initial condition, our initial condition must be $\mathbf{u}^0(x) = (\hat{\mathbf{u}}^0,0)$. We must begin with an initial condition that does not have any unresolved modes activated. For example, we will use $u^0(x) = \sin(x)$, for which only the first Fourier mode is nonzero.

The Markovian term for $\hat{u}_k(t)$ for $k\in F$ is given by $Pe^{t\mathcal{L}}P\mathcal{L}\hat{u}_k^0$:
\begin{align*}
R_k^0(\hat{\mathbf{u}}) = Pe^{t\mathcal{L}}P\mathcal{L}\hat{u}_k^0&=Pe^{t\mathcal{L}}P\left[i\epsilon^2 k^3\hat{u}_k^0-\frac{ik}{2}\sum_{\substack{p+q=k\\p,q\in F\cup G}}\hat{u}_p^0 \hat{u}_q^0\right]\\
R_k^0(\hat{\mathbf{u}})&=i\epsilon^2 k^3\hat{u}_k-\frac{ik}{2}\sum_{\substack{p+q=k\\p,q\in F}}\hat{u}_p \hat{u}_q.
\end{align*}The Markovian term has the same form as the full system, but has been restricted to sums over the resolved modes. We can easily compute the convolution sum in this expression using fast Fourier transforms, but only retaining the resolved modes of the result. In fact, it will be prudent to define a function representing a convolution of two vector valued functions $\mathbf{f}$ and $\mathbf{g}$ with their respective components labeled $f_i$ and $g_i$. We define the convolution of $\mathbf{f}$ with $\mathbf{g}$ with resolved modes retained as:
\begin{equation}
\hat{C}_k(\mathbf{f}(\mathbf{u}),\mathbf{g}(\mathbf{u})) = -\frac{ik}{2}\sum_{\substack{p+q = k\\ k\in F}}f_p(\mathbf{u})g_q(\mathbf{u}).\label{convo}
\end{equation}With this definition, the Markovian term is
\begin{equation}
R_k^0(\hat{\mathbf{u}}) = i\epsilon^2 k^3 \hat{u}_k + \hat{C}_k(\hat{\mathbf{u}},\hat{\mathbf{u}}).\label{kdvMarkov}
\end{equation}For future use, we also define the same convolution, but with only unresolved modes retained. Thus we define:
\begin{equation}\tilde{C}_k(\mathbf{f}(\mathbf{u}),\mathbf{g}(\mathbf{u})) = -\frac{ik}{2}\sum_{\substack{p+q = k\\ k\in G}}f_p(\mathbf{u})g_q(\mathbf{u}).\label{convo_unresolved}
\end{equation}

Finally, it will be useful to define the vector-valued functions $\hat{\mathbf{C}}(\mathbf{f}(\mathbf{u}),\mathbf{g}(\mathbf{u}))$ and $\tilde{\mathbf{C}}(\mathbf{f}(\mathbf{u}),\mathbf{g}(\mathbf{u}))$ whose components are the appropriate convolutions (\ref{convo}) and (\ref{convo_unresolved}). It can be shown that the Markovian term conserves energy in the resolved modes. Thus, it does not allow any transfer of energy out of the resolved modes, which must occur if we are to accurately reproduce what would happen in the full system. That must be accomplished through the memory term.

We can compute terms in the complete memory approximation to generate reduced order models for KdV. We will consider up to fourth order in $t$. This reduced order model for $u_k$ can be written
\begin{equation}
\frac{d\hat{u}_k}{dt} = R_k^0(\hat{\mathbf{u}})+\sum_{i=1}^4\alpha_i(t)t^iR_k^i(\hat{\mathbf{u}}).\label{newROM}
\end{equation}

The Markov term $R_k^0$ was computed above in \eqref{kdvMarkov}. We next compute the $t$-model term $R_k^1(\hat{\mathbf{u}})$. First, we compute $Q\mathcal{L}\hat{u}_k^0$ (again for $k\in F$):
\begin{align*}
Q\mathcal{L}\hat{u}_k^0 &= \mathcal{L}\hat{u}_k^0-P\mathcal{L}\hat{}u_k^0\\
&=-\frac{ik}{2}\sum_{\substack{p+q=k\\p\in F,q\in G}}\hat{u}_p^0 \tilde{u}_q^0-\frac{ik}{2}\sum_{\substack{p+q=k\\p\in G,q\in F}}\tilde{u}_p^0 \hat{u}_q^0-\frac{ik}{2}\sum_{\substack{p+q=k\\p,q\in G}}\tilde{u}_p^0 \tilde{u}_q^0.
\end{align*}Next we compute $P\mathcal{L}Q\mathcal{L}\hat{u}_k^0$:
\begin{align*}
P\mathcal{L}Q\mathcal{L}\hat{u}_k^0
&=-\frac{ik}{2}\sum_{\substack{p+q=k\\p\in F,q\in G}}P\mathcal{L}[\hat{u}_p^0 \tilde{u}_q^0]-\frac{ik}{2}\sum_{\substack{p+q=k\\p\in G,q\in F}}P\mathcal{L}[\tilde{u}_p^0 \hat{u}_q^0]\\
&=2\left(-\frac{ik}{2}\sum_{\substack{p+q=k\\p\in F,q\in G}}\hat{u}_p^0\left(-\frac{iq}{2}\sum_{\substack{r+s=q\\r,s\in F}}\hat{u}_r^0\hat{u}_s^0\right)\right).
\end{align*}

We have here repeatedly used the projection operator $P$ to eliminate terms that become zero. We have also made use of the symmetry of the two sums that remained. Therefore, the $t$-model term is
\begin{equation}
R_k^1(\hat{\mathbf{u}}) = 2\hat{C}_k(\hat{\mathbf{u}},\tilde{\mathbf{C}}(\hat{\mathbf{u}},\hat{\mathbf{u}})).\label{kdvt}
\end{equation}Because the convolutions involve a function with only resolved modes convolved with a function with only unresolved modes, the result is dealiased by construction if we augment our Fourier vectors with one additional mode.

We can compute additional terms entirely in the shorthand established above by recognizing the following rules:
\begin{enumerate}
\item Because a convolution sum is a product of terms and $\mathcal{L}$ is a differential operator, it operates according to the product rule. That is, for every argument in a convolution, we get a term that is a duplicate of that convolution, but with $\mathcal{L}$ applied to that argument. For example:
\begin{align*}
\mathcal{L}\hat{C}_k(\hat{\mathbf{u}}^0,\tilde{\mathbf{C}}(\hat{\mathbf{u}}^0,\hat{\mathbf{u}}^0)) =& \hat{C}_k(\mathcal{L}\hat{\mathbf{u}}^0,\tilde{\mathbf{C}}(\hat{\mathbf{u}}^0,\hat{\mathbf{u}}^0)) \\&+ \hat{C}_k(\hat{\mathbf{u}}^0,\tilde{\mathbf{C}}(\mathcal{L}\hat{\mathbf{u}}^0,\hat{\mathbf{u}}^0)) \\&+ \hat{C}_k(\hat{\mathbf{u}}^0,\tilde{\mathbf{C}}(\hat{\mathbf{u}}^0,\mathcal{L}\hat{\mathbf{u}}^0)).
\end{align*}
\item Each $\mathcal{L}\hat{\mathbf{u}}^0$ term is expanded as \[\mathcal{L}\hat{\mathbf{u}}^0 = i\epsilon^2  \hat{\mathbf{u}}^{k3} + \hat{\mathbf{C}}(\hat{\mathbf{u}}^0,\hat{\mathbf{u}}^0)+2\hat{\mathbf{C}}(\hat{\mathbf{u}}^0,\tilde{\mathbf{u}}^{0})+\hat{\mathbf{C}}(\tilde{\mathbf{u}}^{0},\tilde{\mathbf{u}}^{0}),\] where $(\hat{\mathbf{u}}^{k3})_j = j^3 \hat{u}_j^0$. $\mathcal{L}\tilde{\mathbf{u}}^0$ is expanded in an identical manner, but with each term being the unresolved part, rather than the resolved part.
\item When $\mathcal{L}$ is applied to terms involving powers of $k$, the following occurs:
\begin{align*}
\mathcal{L} (i\epsilon^2\hat{\mathbf{u}}^{k3})_k &= \mathcal{L} i\epsilon^2k^3 \hat{u}_k^0 \\&= (i\epsilon^2 k^3)^2 \hat{u}_k^0 + (i\epsilon^2 k^3)[\hat{C}_k(\hat{\mathbf{u}}^0,\hat{\mathbf{u}}^0)+ 2\hat{C}_k(\hat{\mathbf{u^{0}}},\tilde{\mathbf{u}}^0)+\hat{C}_k(\tilde{\mathbf{u}}^{0},\tilde{\mathbf{u}}^0)]\\&=\left(-\epsilon^4\hat{\mathbf{u}}^{k6}+i\epsilon^2[\hat{\mathbf{C}}^{k3}(\hat{\mathbf{u}}^0,\hat{\mathbf{u}}^0) + 2 \hat{\mathbf{C}}^{k3}(\hat{\mathbf{u}}^0,\tilde{\mathbf{u}}^0) +\hat{\mathbf{C}}^{k3}(\tilde{\mathbf{u}}^0,\tilde{\mathbf{u}}^0)]\right)_k
\end{align*}where $(\hat{\mathbf{u}}^{k6})_j=j^6 \hat{u}_j^0$, $(\hat{\mathbf{C}}^{k3}(\mathbf{a},\mathbf{b}))_j =  j^3\hat{C}_j(\mathbf{a},\mathbf{b})$. That is, the further expansion of powers of $k$ proceed in an easily understandable manner.
\end{enumerate}
The projection operator $P$ is also simple to implement. When applied to a convolution, it applies to each term in the convolution. $P$ sets equal to zero any occurrence of $\tilde{\mathbf{u}}^0$. $Q$ can be represented as $I-P$. Finally, $Pe^{t\mathcal{L}}$ merely advances the initial conditions $\mathbf{u}^0$ to the current time $\mathbf{u}$. We implemented these definitions into a Mathematica notebook, which can then generate ROMs of any order for KdV. In fact, through redefinitions of $P$ and $\mathcal{L}$, this software can be used to derive ROMs of any order for a generic PDE with a generic projector (the software is available upon request).

Using this, we computed expressions for $R_k^2$, $R_k^3$, and $R_k^4$, though they are too large to express here. They can be found in the Appendix. These terms involve convolutions of unresolved terms with other unresolved terms, necessitating that we further augment our vectors using the $3/2$-rule to dealias the results.

The functional form of these higher-order terms in the complete memory approximation are quite complicated, and it is unlikely that a mathematical modeler would propose them. However, because they are derived from the dynamics themselves through the Mori-Zwanzig formalism, we find that they inherit significant structure from the full KdV equation.

\subsection{Renormalization Coefficients}\label{renorm_coeff}

The non-renormalized ROMs as produced directly from the complete memory approximation (\ref{nonrenormalized}) are numerically unstable, so we focus on the renormalized ROMs (\ref{renormalized}). We must develop a procedure for computing the renormalization coefficients $\alpha_i$. An important quantity in a ROM for KdV is the mass in a Fourier mode:
\begin{equation}
M_k(t) =  |u_k(t)|^2.
\end{equation}For select values of $\epsilon$, we computed the exact solution, from which we compute the rate of change of the mass:
\begin{equation}
\Delta M_k(t) = u_k(t)\overline{R_k(\mathbf{u})}+\overline{u_k(t)} R_k(\mathbf{u}).
\end{equation} Contrary to Burgers, the mass in a subset of modes is not monotonically decreasing. Instead, there is a ``mass rebound'' as the dispersive term opposes the formation of a shock, and mass returns from high-frequency modes to low-frequency modes. 

We also computed $R_k^i(\mathbf{u})$ for ROMs constructed for those subsets. From this, we can measure the impact each $R_k^i$ term has upon the rate of change of the mass in individual modes. This is given by:

\begin{equation}
\Delta M_k^i(t) = \hat{u}_k(t) \overline{R_k^i(\hat{\mathbf{u}})}+\overline{\hat{u}_k(t)} R_k^i(\hat{\mathbf{u}}).
\end{equation}

Any net rate of change of mass in the resolved modes must be accounted for by these memory terms alone, because the Markovian term conserves mass in the resolved modes.

\begin{figure}
\centering
\includegraphics[width=\linewidth]{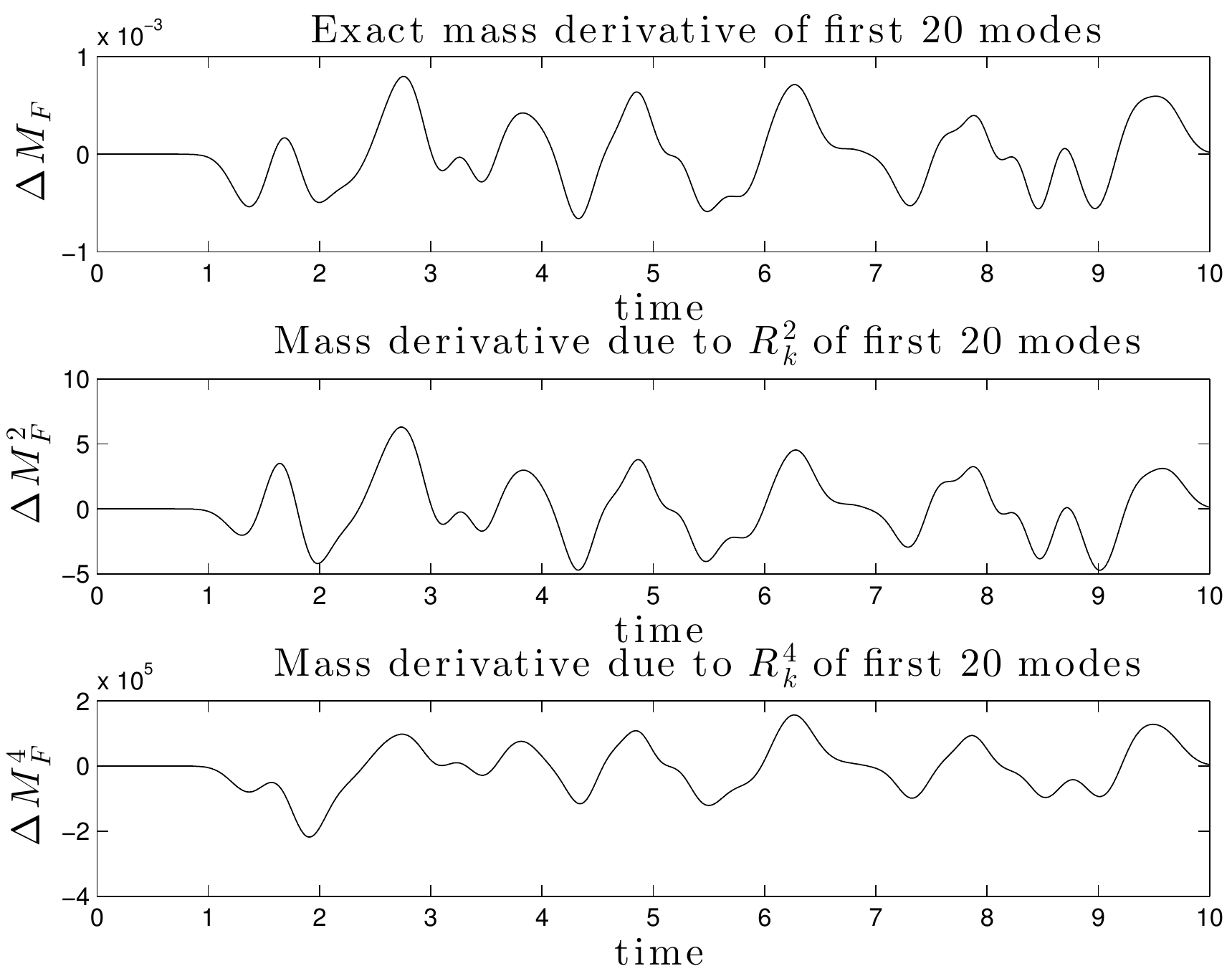}
\caption{(Top) The exact net mass derivative of the first $N=20$ positive modes of the solution to the KdV equation with $\epsilon = 0.1$ up to time 10. (Middle) The rate of change of the mass in the first $N=20$ positive modes for the $R_k^2$ term in an ROM of size $N=20$ provided the exact solution $\hat{\mathbf{u}}$ as input. (Bottom) The rate of change of the mass in the first $N=20$ positive modes for the $R_k^4$ term in an ROM of size $N=20$ provided the exact solution $\hat{\mathbf{u}}$ as input. Note that the ROM curves seem to differ from the exact solution by only a constant scaling factor.}
\label{fig:energy_compare}
\end{figure}

We found that the $\Delta M_k^2$ and $\Delta M_k^4$ terms closely mirror the exact net mass derivative $\Delta M_k$, as depicted in Figure \ref{fig:energy_compare}.

From this, we draw several conclusions. First, this close agreement suggests that the complete memory approximation is in some sense a ``correct'' way of expanding the memory integral. Second, because the curves appear to differ by only a constant factor, it suggests a functional form for the renormalization coefficients of
\begin{equation}
\alpha_i(t) = \alpha_i t^{-i}, \label{coeff}
\end{equation}
such that the time dependence of the renormalization coefficient cancels the time dependence in the memory terms. Note that the coefficients $\alpha_i(t)$ exhibit algebraic temporal dependence which characterizes the {\it absence} of timescale separation between the resolved and unresolved variables. Also, that the reduced-order models we employ do not require an integral memory convolution, yet they incorporate long memory effects.

For a fixed $\epsilon$ and simulation length $T$, we used a least squares fit to identify optimal choices for the constants $\alpha_i$. Suppose we have the exact solution at times $\{t_0,\dots,t_m\}$. The cost function for a given set of resolved modes $F$ is:

\begin{align}
C_F(\mathbf{\alpha}) =& \sum_{j=0}^m\sum_{k\in F} \left[\Delta M_k(t_j) - \sum_{i=1}^4 \alpha_i\Delta M_k^i(t_j)\right]^2\notag\\
+&\sum_{j=0}^m\left[\sum_{k\in F} \biggl( \Delta M_k(t_j) - \sum_{i=1}^4 \alpha_i\Delta M_k^i(t_j)\biggr) \right]^2.\label{least_sq}
\end{align}That is, we sought to minimize the error in representing the derivative of the mass of each individual mode, as well as the net mass flow in and out of the resolved modes at all times.

We conducted least squares fits for many sets of resolved modes and choices of dispersion $\epsilon$. In each case, we used data for $t_i\in[t_a,t_b]$ in increments of $0.001$. In practice, we found the results were insensitive to the details of the data used. As long as the width of the interval, $t_b-t_a$, was greater than 3 time units, the coefficients were not substantially different. For our calculations here, we used $t_a=0$ and $t_b=10$. We found that the optimal coefficients displayed power law behavior {\it both} in $\epsilon$ and in the size of the reduced order model $N.$ We used an additional least squares fit to identify the scaling coefficients for each variable. 

The scaling law behavior of the coefficients suggests the presence of incomplete similarity for appropriate non-dimensional parameters which we now describe \cite{barenblatt2003}. We begin with the non-dimensionalization of the KdV equation \eqref{kdv_equation}. If $L$ and $U$ are characteristic length and velocity scales then we can define the non-dimensional variables $x'=x/L,$ $u'=u/U,$ and $t'=t\frac{U}{L}=t/T.$ We chose the length of the domain for the characteristic length $L$ while for the characteristic velocity we used $U=[\frac{\int_0^L u^2_0(x)dx}L]^{1/2}.$ After dropping the primes, Eq. \eqref{kdv_equation} becomes
\begin{equation}
u_t + uu_x +\frac{1}{Re^2}u_{xxx}=0, \label{kdv_equation_non}
\end{equation}
where $Re=\frac{\sqrt{U}L}{\epsilon}$ is a ``dispersive'' Reynolds number. From Eq. \eqref{renormalized} we know that the renormalized coefficients $\alpha_i(t)$ are non-dimensional. This is because the terms $t^i R^i_k$ that they multiply in \eqref{renormalized} have the correct dimensions by construction, since they are produced through the MZ formalism and not added by hand. Also, from Eq. \eqref{coeff} we know that the prefactors $\alpha_i$ have dimension $T^i=(\frac{L}{U})^i.$ If in addition to the Reynolds number $Re$ we define the non-dimensional parameters $\Pi_i=\frac{\alpha_i}{(\frac{L}{U})^i}$ and $\Lambda=\frac{N}{L^{-1}}$ we can represent the renormalized coefficients as
\begin{equation}
\Pi_i = a_i Re^{b_i}\Lambda^{c_i}. \label{t4_renorm}
\end{equation}

\begin{figure}
\begin{center}
\includegraphics[width = 0.8\textwidth]{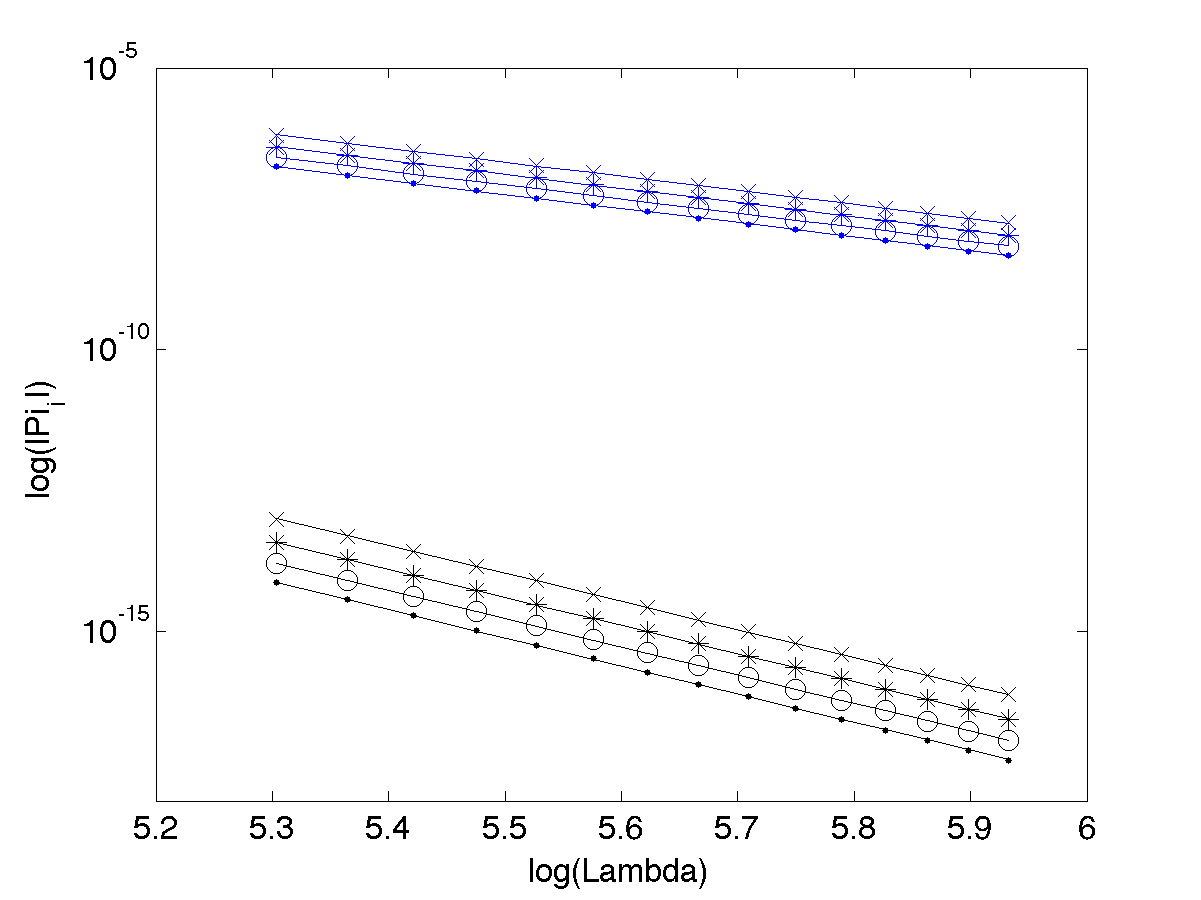}

\includegraphics[width = 0.8\textwidth]{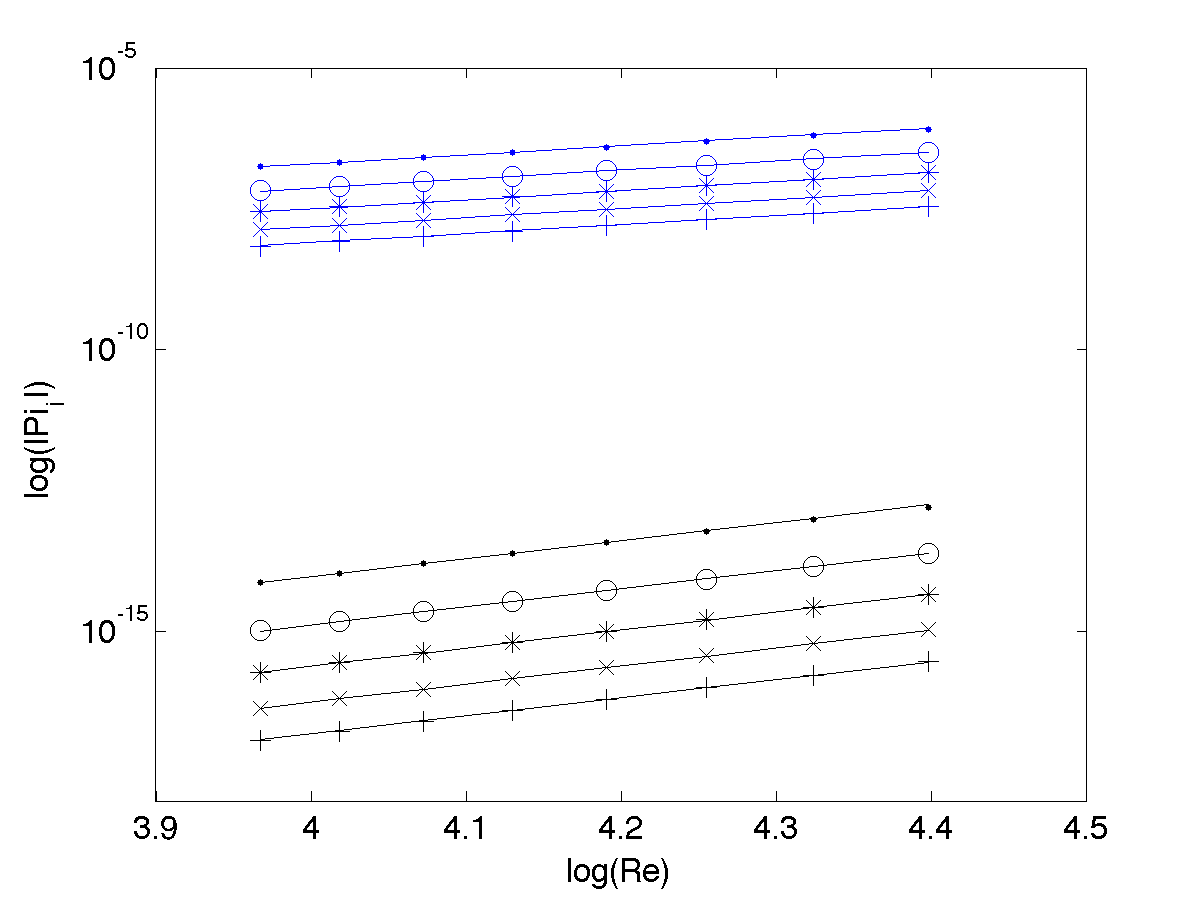}
\end{center}
\caption{Log-log plots of the absolute values of the nondimensionalized renormalization coefficients $\Pi_i$. Computed by minimizing (\ref{least_sq}) using data from time $[0,10]$ for a variety of values of dispersive Reynolds numbers $Re$ and resolutions $\Lambda$ with an initial condition $u_0(x)=\sin(x)$. (Top) $\log(\Lambda)$ plotted against $\log(|\Pi_2|)$ (blue) and $\log(|\Pi_4|)$ (black) for Reynolds numbers corresponding to $\epsilon$ equal to 0.1 (dots), 0.09 (circles), 0.08 (stars), and 0.07 (crosses). (Bottom) $\log(Re)$ plotted against $\log(|\Pi_2|)$ (blue) and $\log(|\Pi_4|)$ (black) for $\Lambda$ corresponding to $N$ equal to 32 (dots), 38 (circles), 44 (stars), 50 (crosses), and 56 (pluses).}\label{fig:power_laws}
\end{figure}

Figure \ref{fig:power_laws} presents the results of the least squares fit in appropriate non-dimensional variables. For $i=1$ and $i=3$, we found $a_i=0$, so the power law exponents $b_i$ and $c_i$ were irrelevant. It is notable that the odd numbered terms in the memory expansion seem to be unnecessary for capturing the dynamics (we will discuss this more in Section \ref{discussion}). For $i=2$, we found $a_2 = -1.2473$, $b_2 = 3.6910$, and $c_2 = -5.7356$. For $i=4$, we found $a_4 =  -0.3675$, $b_4 = 7.3881$, and $c_4 = -11.4719$ (we will refer to reduced order models of this form with these renormalization coefficients the fourth order complete renormalized model). We note that $b_4$ and $c_4$ are roughly double $b_2$ and $c_2$ respectively. It seems likely that these scaling laws belie an even deeper structure in the memory inherited by the ROMs from the original KdV equation than we have uncovered here. This is a subject of future investigation by the authors.

The non-dimensional renormalization constants $\Pi_i$ in \eqref{t4_renorm} exhibit incomplete similarity in {\it both} parameters $Re$ and $\Lambda.$ If we keep the resolution of the reduced model fixed ($\Lambda$ fixed), then $\Pi_i \rightarrow \infty$ as $Re \rightarrow \infty.$ This is a hallmark of a singular perturbation problem which is to be expected since for $Re=\infty$ we recover the Burgers equation which develops shocks in finite time, while KdV does not. Similarly, if we keep $Re$ fixed, and let the resolution of the reduced model $\Lambda \rightarrow \infty$ then $\Pi_i \rightarrow 0$ which is also expected since in this limit there is no need for a memory term.

We used a similar approach to find renormalization coefficients $\beta_i$ for a ROM  of the form
\[\frac{d\hat{u}_k}{dt} = R_k^0(\hat{\mathbf{u}})+\sum_{i=1}^2 \beta_i(t) t^i R_k^i(\hat{\mathbf{u}}).
\]Similarly, the time dependence of the coefficients was $\beta_i(t) = \beta_i t^{-i}$. Let $\Pi_i'$ be the non-dimensionalized form of $\beta_i$. It also obeys a scaling law formula:
\begin{equation}
\Pi'_i = d_i Re^{e_i}\Lambda_1^{f_i},\label{t2_renorm}
\end{equation}where $d_1 =0$ (so the $t$-model again is not included), and $d_2 = -0.7615$, $e_2 = 3.7681$, and $f_2 = -5.8081$. Observe that the scaling coefficients for $\Lambda_1$ and $Re$ are very similar to the results for the second order term in the fourth order model, but that the prefactor is different (we will call reduced order models with this form and renormalization coefficients the second order complete renormalized model). Once again, the fact that the $t$-model is unimportant for capturing the memory effects of the KdV equation is noteworthy, as is the scaling structure inherent in the $t^2$-model coefficient.

Finally, we should note that the renormalization constant for the $R_k^2$ term is negative in both cases, though the non-renormalized $R_k^2$ term closely mirrors the rate of change of the energy in the resolved modes. The sign change causes this agreement to be lost, though the change of energy in each individual mode is well-captured. This curious phenomenon is also the subject of further investigation.

\section{Results}

Because a fully resolved solution is possible for sufficiently large values of $\epsilon$, we are able to test these new renormalized reduced order models for both stability and accuracy by comparing them against the exact solution. A natural comparison for our ROMs is the ``average dynamics'' given by the Markov term alone, disregarding memory terms. All simulations are conducted using a time step of $\Delta t = 0.001$ using an implicit-explicit integration scheme that has proven effective for the KdV equation with periodic boundary conditions \cite{klein2008fourth,driscoll2002composite}.

There are several metrics for identifying the accuracy of a reduced order model of this kind. First, we have seen from the exact solution that for a given subset of modes $F$ in the exact solution, mass should flow both in and out as time passes. This corresponds to mass passing from the resolved modes to the unresolved modes and back. A ROM for the same subset of modes should accurately capture this net in and out flow of mass. The Markov term alone conserves mass in the resolved modes, and so is unable capture this effect at all. 

\begin{figure}[h]
\begin{center}
\includegraphics[width=0.8\textwidth]{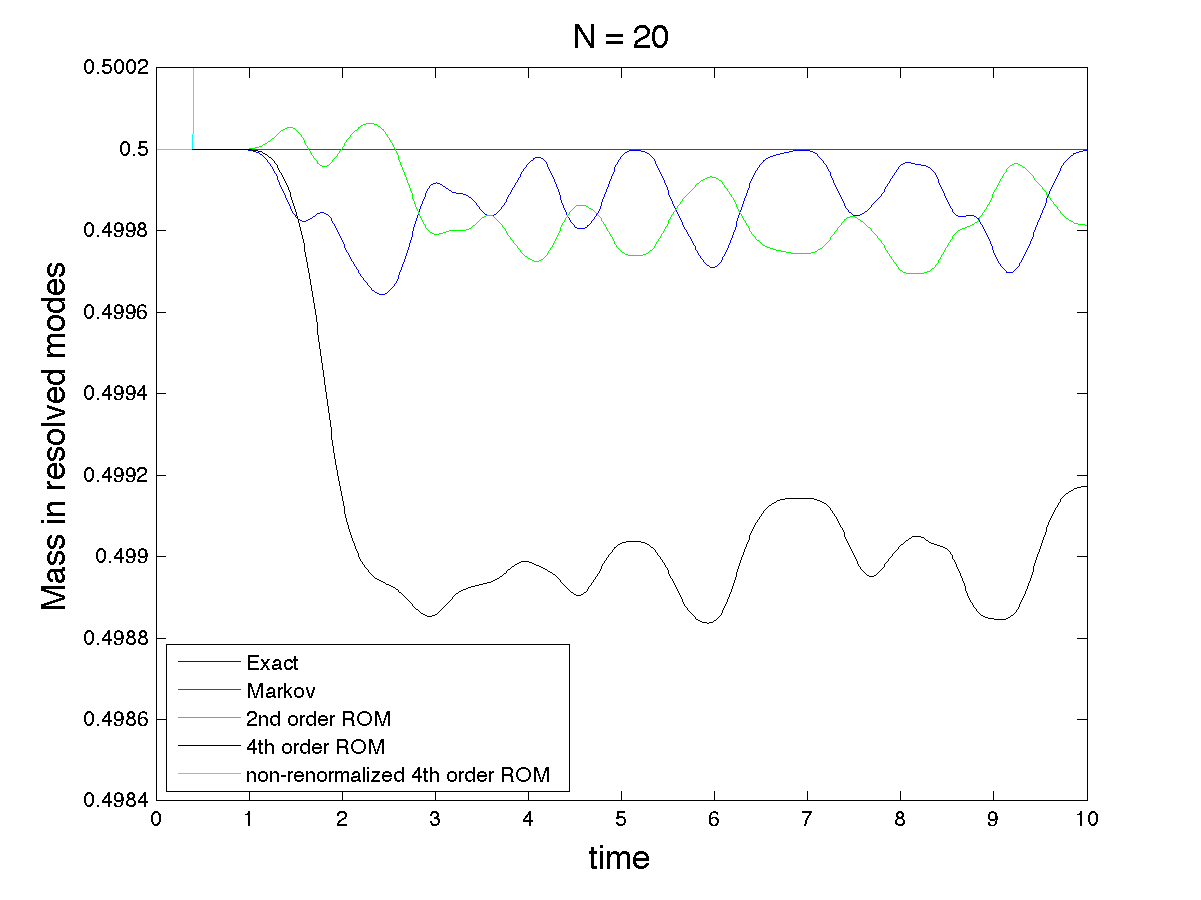}
\end{center}
\caption{The mass in the first $N=20$ positive fourier modes for an initial condition of $\sin(x)$ according to several models. The blue curve depicts the exact solution, found by running a simulation with $N = 256$ positive modes. The red curve depicts the Markov model with $N=20$ modes. The green curve depicts the renormalized second order complete ROM, while the black curve depicts the renormalized fourth order complete memory ROM. The cyan curve is a non-renormalized version of the 4th order ROM, which is unstable. }\label{fig:energy_over_time}
\end{figure}

An example of this dynamic mass loss and gain is depicted in Figure \ref{fig:energy_over_time}. Observe that the second order ROM drains mass when the exact system should be gaining it, and vice versa. This is a consequence of the fact that the renormalization constant is negative, while the $R_k^2$ term itself appears to be a scaled version of the derivative of the total mass in the resolved modes (as seen in Figure \ref{fig:energy_compare}). The fourth order ROM, on the other hand, drains too much mass in the initial three units of time, but then gains and drains mass consistently with the exact rate. The addition of the fourth order term seems to improve the model's ability to capture mass transfer in and out of the resolved modes. Note, finally, how little the mass in the resolved modes actually changes. In the $N=20$ case, less than $0.08\%$ ever leaves the resolved modes. For this reason, intuition would suggest that the memory is not important, and that the Markov model will be sufficient. This proves to not be the case. Also, it is remarkable that this little amount of mass flow in and out of the resolved variables carries with it the rich structure for the renormalized coefficients presented in Section \ref{renorm_coeff}.

\begin{figure}[h]
\begin{center}
\includegraphics[width=0.8\textwidth]{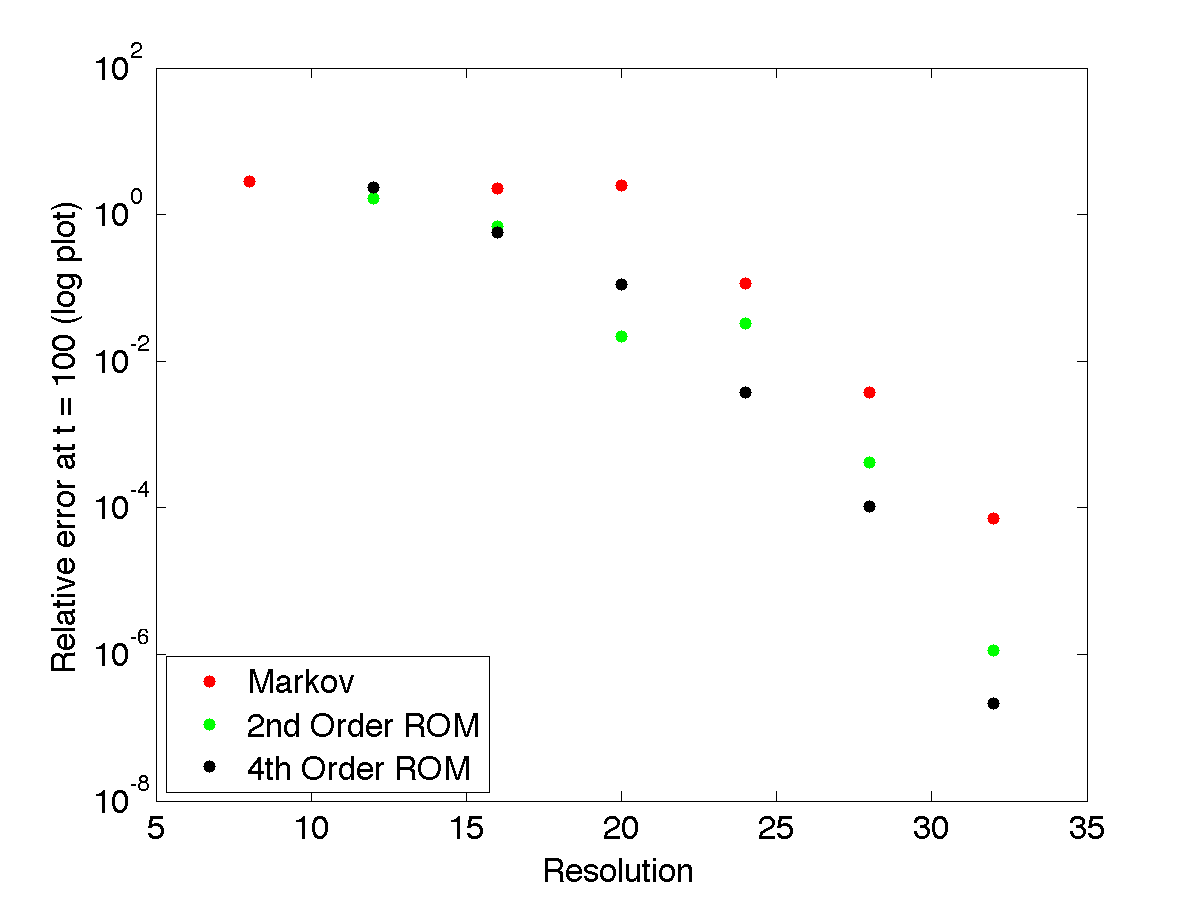}
\end{center}
\caption{The relative global error at time $t=100$ for several models plotted on a logarithmic axis with $\epsilon=0.1$ and initial condition $\sin(x)$. The error is computed as the ratio between the $L_2$ norm of the real space error and the $L_2$ norm of the exact real space solution. All three ROMs improve as the number of resolved modes increases, but the complete memory ROMs are far more accurate than the Markov model alone.  }\label{fig:real_space_error}
\end{figure}

The need for the inclusion of memory can be seen by inspecting our solutions in real space. We compute the $L_2$ norm of the difference between the exact real space solution and that predicted by the ROM. We divide this error by the $L_2$ norm of the exact solution, producing a global relative error between the exact and approximated solution trajectories. Here, it is reasonable to also compute the error of the Markov approximation for comparison. Figure \ref{fig:real_space_error} depicts the global relative error at time $t=100$ for ROMs of several different resolutions with $\epsilon = 0.1$ (the error results for other values of $\epsilon$ are qualitatively similar). Recall that the formulae for the renormalization coefficients were fit only with data between $t\in[0,10]$, so it appears that these coefficients are valid for a long time. 

The accuracy of complete ROMs is not achieved until $N$ is sufficiently large. For $\epsilon=0.1$ this is around $N=20$. For $\epsilon = 0.09$, it is approximately $N=24$. Qualitatively, our results suggest that a stable and accurate ROM can only be constructed for sets of resolved modes whose ``full models'' comprise at least half of the modes needed for a fully resolved simulation. Because KdV is a dispersive problem, the total mass is conserved. Therefore, a reduced order model can only hope to capture the dynamic flow of a fixed amount of mass if the full model upon which it is based contains a large proportion of the modes through which the mass must flow.

As $\epsilon\to0$, the initial draining of mass passes to higher frequencies before beginning to rebound. When $\epsilon=0$, the problem becomes Burgers' equation, which has a finite time singularity, and the mass (in the Burgers literature called energy) cascades to higher frequencies indefinitely. For this case, a renormalized model in which the coefficients $\alpha_i(t)$ are actually {\it constant} performs very well  \cite{stinis2013renormalized}. As $\epsilon\to0$, the renormalization coefficients as computed by our formulae (\ref{t4_renorm}) and (\ref{t2_renorm}) grow to infinity. This suggests that, as $\epsilon$ becomes small, the assumptions made in deriving our complete renormalized ROMs fail to hold. As mentioned before, this is an indication of the singular nature of the perturbation problem as $\epsilon \to 0.$

Finally, we must comment on the computational cost of these reduced order models. As observed, stable and accurate ROMs can only be constructed for resolutions that are approximately four times smaller than a fully resolved solution. The FFTs employed in an ROM are large enough to contain the ``full model,'' which is twice as large as the ROM itself. Thus, in the best case scenario, we will be using FFTs that are half as large as those needed for the full solution. The second order ROM uses six FFTs and IFFTs per timestep, while the fourth order ROM requires 22. A fully resolved simulation, on the other hand only uses one FFT and IFFT per timestep, though it is twice as large. Consequently, the ROMs produced in this paper are necessarily less efficient than simply solving the problem outright. Their value is not found in producing efficient computational schemes, but in shedding light on the role of memory in dispersive problems.

\section{Discussion and future work}\label{discussion}

We have developed a new method of constructing reduced order models from the Mori-Zwanzig formalism. This method involves expanding the memory term in a Taylor series in time, and repeatedly applying the Mori-Zwanzig formalism to unprojected terms. Different ROMs can be constructed by truncating the Taylor series at different terms. The telescoping terms in the expansion are challenging to derive, but the process can be automated in a symbolic software package. The functional forms of these terms are complex, and would be unlikely to be discovered by classical mathematical modeling techniques. Instead, these terms are derived from the dynamics themselves.

The resulting ROMs are numerically unstable, unless each term is tamed with a specific renormalization coefficients. We find that these renormalization coefficients decay algebraically in time, canceling the time dependence originally found in the Taylor expansions. This differs from previous efforts to renormalize ROMs derived from the Mori-Zwanzig formalism, which became stable with renormalization coefficients that were constant in time.

The renormalization constants that most effectively capture the correct dynamics of the resolved modes also suggest significant structure in the memory. The odd terms in the memory expansion seem to be unimportant for capturing the dynamics of the memory. This is at odds with previous results for dissipative problems, such as Burgers' equation, for which the odd terms seemed to be more important \cite{stinis2015renormalized}. It is possible that the even terms in the memory expansion account for the memory of dispersive terms, while the odd terms account for the memory of dissipative terms. To clarify the situation, the authors plan to apply it to dissipative problems (such as Burgers' equation) and problems with dissipation and dispersion (such as the KdV-Burgers equation). For the case of KdV-Burgers equation, averaging and renormalization for traveling waves was considered in \cite{chorin2003kdvburgers}. 

The structure of the memory is also evident in the fact that the coefficients of the even terms in the memory expansion obey strict scaling laws in $\epsilon$ and $N$, with the power law dependence for the $R_k^4$ term being approximately twice as large as that of the $R_k^2$ term. Due to this structure, we can construct ROMs for any $\epsilon$ and $N$ without needing to know the exact solution for that specific $\epsilon$. The highly structured nature of the renormalization coefficients suggest that the complete memory approximation is, in some manner, a natural way of describing the dynamics of a subset of modes in this problem.

We found qualitatively that accurate ROMs were only possible when the ``full model'' on which the reduced model is based is nearly fully resolved itself. This was not the case for dissipative problems, and we hypothesize that it is the dispersive nature of the problem that necessitates large ``full models.'' Once accurate ROMs can be constructed, however, they prove to be much more effective than the Markov term alone. This is noteworthy, because the amount of mass leaving the resolved modes is exceptionally small. This indicates that, even when the memory term is extremely small, ignoring it will lead to inaccurate results, especially during long time simulations. This is a hallmark of a singular perturbation problem (in our case the perturbation is the small magnitude of the dispersion coefficient $\epsilon$). Also, this agrees with the discovered incomplete similarity exhibited by the renormalized coefficients \cite{barenblatt2003}.

The results of the complete memory approximation and its renormalization are promising and important. First, they support the need for memory terms in reduced models. Second, the suggested distinction of different order MZ memory terms as relevant for different physical mechanisms points towards a systematic classification of the memory. This realization supports the advantage of model reduction starting from an exact formalism like MZ over adding terms by hand. Finally, the success of renormalization in producing stable and accurate models for long times indicates an elegant and effective way of incorporating dynamic information about the scales we aim to resolve.

\section{Acknowledgements}

JP would like to thank Prof. B. Deconinck for suggesting the KdV equation as a test case for renormalized reduced order modeling, and B. Segal and L. Thompson for helpful and interesting discussions. PS would like to thank Profs. K. Lin, J. Harlim, F. Lu and K. Duraisamy for very useful discussions. The work of PS was partially supported by the U.S. Department of Energy Office of Science, Office of Advanced Scientific Computing Research, Applied Mathematics program, Collaboratory on Mathematics for Mesoscopic Modeling of Materials (CM4), under Award Number DE-SC0009280. 

\appendix

\section{Complete Memory Approximation}

Once the evolution operators in the memory integral of the Mori-Zwanzig formalism (\ref{reduced_MZ}) are expanded and integrated termwise, we are left with (\ref{eq:novel}), which has been reproduced here:

\begin{equation}P\int_0^te^{(t-s)\mathcal{L}}P\mathcal{L}e^{sQ\mathcal{L}}Q\mathcal{L}u_k^0\,\mathrm{d}s\notag\\
=Pe^{t\mathcal{L}}\left(\sum_{i=0}^\infty \sum_{j=0}^\infty \frac{(-1)^it^{i+j+1}}{i!j!(i+j+1)}\mathcal{L}^iP\mathcal{L}(Q\mathcal{L})^{j+1}u_k^0\right).\label{appendix_complete}
\end{equation}The $O(t)$ term is the $t$-model $tPe^{t\mathcal{L}}[P\mathcal{L}Q\mathcal{L}]u_k^0$, and the $O(t^2)$ term is
\[\frac{t^2}{2}Pe^{t\mathcal{L}}[P\mathcal{L}Q\mathcal{L}Q\mathcal{L}-\mathcal{L}P\mathcal{L}Q\mathcal{L}]u_k^0.
\]In order to produce a term that is projected \emph{prior} to evolution, we expand the second term in the parentheses according to its own Mori-Zwanzig formalism, with its memory term also computed as a termwise integrated double sum. This yields
\begin{align}
Pe^{t\mathcal{L}}&\mathcal{L}P\mathcal{L}Q\mathcal{L}u_k^0  = \frac{\partial}{\partial t}Pe^{t\mathcal{L}}P\mathcal{L}Q\mathcal{L}u_k^0\notag\\
&= Pe^{t\mathcal{L}}P\mathcal{L}P\mathcal{L}Q\mathcal{L}u_k^0 + Pe^{t\mathcal{L}}\left(\sum_{i=0}^\infty\sum_{j=0}^\infty \frac{(-1)^it^{i+j+1}}{i!j!(i+j+1)}\mathcal{L}^iP\mathcal{L}(Q\mathcal{L})^{j+1}P\mathcal{L}Q\mathcal{L}u_k^0\right).\label{t2_problem_expanded}
\end{align}Recall that this term is being multiplied by $t^2/2$, so only the first term in (\ref{t2_problem_expanded}) remains $O(t^2)$. Thus, the $O(t^2)$ of the  complete memory approximation is
\begin{equation}
-\frac{t^2}{2}Pe^{t\mathcal{L}}[P\mathcal{L}P\mathcal{L}Q\mathcal{L}-P\mathcal{L}Q\mathcal{L}Q\mathcal{L}]u_k^0.\label{t2}
\end{equation}In order to account for $O(t^3)$ we must also include the $O(t)$ term of (\ref{t2_problem_expanded}), since it becomes $O(t^3)$ when multiplied by $t^2/2$. In addition, some of the $O(t^3)$ from the expansion of the sum (\ref{appendix_complete}) have a leading $\mathcal{L}$, and must be expanded in its own memory expansion in the same manner as above. Once this is done, the $O(t^3)$ term is:
\begin{equation}
\frac{t^3}{6}[P\mathcal{L}P\mathcal{L}P\mathcal{L}Q\mathcal{L}-2P\mathcal{L}P\mathcal{L}Q\mathcal{L}Q\mathcal{L}-2P\mathcal{L}Q\mathcal{L}P\mathcal{L}Q\mathcal{L}+P\mathcal{L}Q\mathcal{L}Q\mathcal{L}Q\mathcal{L}]u_k^0.\label{t3}
\end{equation}Finally, we can repeat these same steps to find the $O(t^4)$:

\begin{align}
-\frac{t^4}{24}\bigg[&P\mathcal{L}Q\mathcal{L}Q\mathcal{L}Q\mathcal{L}Q\mathcal{L}-3P\mathcal{L}Q\mathcal{L}Q\mathcal{L}P\mathcal{L}Q\mathcal{L}-5P\mathcal{L}Q\mathcal{L}P\mathcal{L}Q\mathcal{L}Q\mathcal{L}\notag\\
&-3P\mathcal{L}P\mathcal{L}Q\mathcal{L}Q\mathcal{L}Q\mathcal{L}+3P\mathcal{L}P\mathcal{L}P\mathcal{L}Q\mathcal{L}Q\mathcal{L}+5P\mathcal{L}P\mathcal{L}Q\mathcal{L}P\mathcal{L}Q\mathcal{L}\notag\\
&+3P\mathcal{L}Q\mathcal{L}P\mathcal{L}P\mathcal{L}Q\mathcal{L}-P\mathcal{L}P\mathcal{L}P\mathcal{L}P\mathcal{L}Q\mathcal{L}\bigg]u_k^0.\label{t4}
\end{align}

\subsection{The BCH Approximation}\label{bch_expansion}

In previous work \cite{stinis2013renormalized}, a different method was used to approximate the memory integral. We rewrite the memory term by reversing our use of Dyson's formula
\begin{align*}
P\int_0^te^{(t-s)\mathcal{L}}P\mathcal{L}e^{sQ\mathcal{L}}Q\mathcal{L}u_k^0\,\mathrm{d}s=&Pe^{t\mathcal{L}}Q\mathcal{L}u_k^0-Pe^{tQ\mathcal{L}}Q\mathcal{L}u_k^0\\
=&Pe^{t\mathcal{L}}\left(Q\mathcal{L}u_k^0-e^{-t\mathcal{L}}e^{tQ\mathcal{L}}QLu_k^0\right)\\
=&Pe^{t\mathcal{L}}\left(Q\mathcal{L}u_0^k -e^{C(t)}Q\mathcal{L}u_0^k\right)
\end{align*}where $C(t) = -tP\mathcal{L}+[tP\mathcal{L},tQ\mathcal{L}]+\dots$ is the BCH series, which expands in powers of the commutator $[tP\mathcal{L},tQ\mathcal{L}]$. In the event that $[tP\mathcal{L},tQ\mathcal{L}]$ is small, which corresponds to cases in which the orthogonal and the projected dynamics approximately commute, $C(t)$ can be approximated by $-tP\mathcal{L}$. This may be the case for particular initial conditions. The resulting approximation is:
\begin{equation}
P\int_0^te^{(t-s)\mathcal{L}}P\mathcal{L}e^{sQ\mathcal{L}}Q\mathcal{L}u_k^0\,\mathrm{d}s\approx \sum_{j=1}^\infty (-1)^{j+1}\frac{t^j}{j!}Pe^{t\mathcal{L}}(P\mathcal{L})^jQ\mathcal{L}u_k^0.\label{BCHROM}
\end{equation}

We designate this method of approximating the memory term the ``BCH approximation.'' If we allow $P\mathcal{L}$ and $Q\mathcal{L}$ to commute, we can arrange the terms derived from the complete memory approximation (\ref{t2}), (\ref{t3}), and (\ref{t4}) into the form expressed in (\ref{BCHROM}). Thus, the BCH approximation represents a special case of the complete memory approximation.

\section{Complete Memory Approximation of KdV}

In order to simulate ROMs for KdV constructed through the complete memory approximation, we must use the definition of $\mathcal{L}$ and $P$ to construct the $t$-model term, and the higher order terms (\ref{t2}), (\ref{t3}), and (\ref{t4}). We derived the form of the $t$-model term (\ref{kdvt}) and a set of rules for applying $\mathcal{L}$ and $P$ to the convolution sums that arise in these expansions.

Section \ref{kdv} describes these details. The process becomes quite exhausting, because each nested convolution expands into multiple terms with additional applications of $\mathcal{L}$. Each term in the complete memory approximation also contains more terms than those that preceded it. The process of deriving and simplifying these terms has been automated in a Mathematica notebook.

We will demonstrate the derivation of the $t^2$-model term for a resolved mode $\hat{u}_k$ here as an example. We have already computed $Q\mathcal{L}\hat{u}_k^0$, but future applications of $Q\mathcal{L}$ will be computed by $Q\mathcal{L}=\mathcal{L}-P\mathcal{L}$. Thus, the $t^2$-model term is:
\[Pe^{t\mathcal{L}}[P\mathcal{L}P\mathcal{L}Q\mathcal{L}-P\mathcal{L}Q\mathcal{L}Q\mathcal{L}]\hat{u}_k^0 = Pe^{t\mathcal{L}}[2P\mathcal{L}P\mathcal{L}Q\mathcal{L}-P\mathcal{L}\mathcal{L}Q\mathcal{L}]\hat{u}_k^0.
\]The first term can be computed by applying $P\mathcal{L}$ to $P\mathcal{L}Q\mathcal{L}\hat{u}_k^0$, which we have already found:
\begin{align*}
P\mathcal{L}P\mathcal{L}Q\mathcal{L}\hat{u}_k^0 =& P\mathcal{L}\left[2\hat{C}_k(\hat{\mathbf{u}}^0,\tilde{\mathbf{C}}(\hat{\mathbf{u}}^0,\hat{\mathbf{u}}^0))\right]\\
=&P\bigg[2\hat{C}_k(i\epsilon^2 \hat{\mathbf{u}}^{k3} + \hat{\mathbf{C}}(\hat{\mathbf{u}}^0,\hat{\mathbf{u}}^0)+2\hat{\mathbf{C}}(\hat{\mathbf{u}}^0,\tilde{\mathbf{u}}^0)+\hat{\mathbf{C}}(\tilde{\mathbf{u}}^0,\tilde{\mathbf{u}}^0),\tilde{\mathbf{C}}(\hat{\mathbf{u}}^0,\hat{\mathbf{u}}^0))\\
&\quad \;+4\hat{C}_k(\hat{\mathbf{u}}^0,\tilde{\mathbf{C}}(i\epsilon^2 \hat{\mathbf{u}}^{k3} + \mathbf{C}(\hat{\mathbf{u}}^0,\hat{\mathbf{u}}^0)+2\hat{\mathbf{C}}(\hat{\mathbf{u}}^0,\tilde{\mathbf{u}}^0)+\hat{\mathbf{C}}(\tilde{\mathbf{u}}^0,\tilde{\mathbf{u}}^0),\hat{\mathbf{u}}^0))\bigg]\\
=&2\hat{C}_k(i\epsilon^2 \hat{\mathbf{u}}^{k3}+\mathbf{C}(\hat{\mathbf{u}}^0,\hat{\mathbf{u}}^0),\mathbf{C}^*(\hat{\mathbf{u}}^0,\hat{\mathbf{u}}^0))+4\hat{C}_k(\hat{\mathbf{u}}^0,\tilde{\mathbf{C}}(i\epsilon^2 \hat{\mathbf{u}}^{k3}+\hat{\mathbf{C}}(\hat{\mathbf{u}}^0,\hat{\mathbf{u}}^0),\hat{\mathbf{u}}^0)).
\end{align*}The second term in the $O(t^2)$ expansion is:
\begin{align*}
P\mathcal{L}\mathcal{L}Q\mathcal{L}\hat{u}_k^0 =&P\mathcal{L}\mathcal{L}[2\hat{C}_k(\hat{\mathbf{u}}^0,\tilde{\mathbf{u}}^0)+\hat{C}_k(\tilde{\mathbf{u}}^0,\tilde{\mathbf{u}}^0)]\\
=&P\mathcal{L}\bigg[2\hat{C}_k(i\epsilon^2 \hat{\mathbf{u}}^{k3} + \hat{\mathbf{C}}(\hat{\mathbf{u}}^0,\hat{\mathbf{u}}^0)+2\hat{\mathbf{C}}(\hat{\mathbf{u}}^0,\tilde{\mathbf{u}}^0)+\hat{\mathbf{C}}(\tilde{\mathbf{u}}^0,\tilde{\mathbf{u}}^0),\tilde{\mathbf{u}}^0)\\
&\quad \;\;\; + 2\hat{C}_k(\hat{\mathbf{u}}^0,i\epsilon^2 \tilde{u}^{k3} + \tilde{\mathbf{C}}(\hat{\mathbf{u}}^0,\hat{\mathbf{u}}^0)+2\tilde{\mathbf{C}}(\hat{\mathbf{u}}^0,\tilde{\mathbf{u}}^0)+\tilde{\mathbf{C}}(\tilde{\mathbf{u}}^0,\tilde{\mathbf{u}}^0))\\
&\quad \;\;\; + 2\hat{C}_k(\tilde{\mathbf{u}}^0,i\epsilon^2 \hat{\mathbf{u}}^{*k3} + \tilde{\mathbf{C}}(\hat{\mathbf{u}}^0,\hat{\mathbf{u}}^0)+2\tilde{\mathbf{C}}(\hat{\mathbf{u}}^0,\tilde{\mathbf{u}}^0)+\tilde{\mathbf{C}}(\tilde{\mathbf{u}}^0,\tilde{\mathbf{u}}^0))
\bigg].
\end{align*}We will apply the projection $P$ in the same step that we apply $\mathcal{L}$ to save space writing out terms that will be eliminated. The result is:
\begin{align*}
P\mathcal{L}\mathcal{L}Q\mathcal{L}\hat{u}_k^0 =& 4\hat{C}_k(i\epsilon^2 \hat{\mathbf{u}}^{k3} + \hat{\mathbf{C}}(\hat{\mathbf{u}}^0,\hat{\mathbf{u}}^0),\tilde{\mathbf{C}}(\hat{\mathbf{u}}^0,\hat{\mathbf{u}}^0))\\
&+ 2\hat{C}_k(\hat{\mathbf{u}}^0,i\epsilon^2 \tilde{\mathbf{C}}^{k3}(\hat{\mathbf{u}}^0,\hat{\mathbf{u}}^0)+2\tilde{\mathbf{C}}(\hat{\mathbf{u}}^0,i\epsilon^2 \hat{\mathbf{u}}^{k3}+\hat{\mathbf{C}}(\hat{\mathbf{u}}^0,\hat{\mathbf{u}}^0)+\tilde{\mathbf{C}}(\hat{\mathbf{u}}^0,\hat{\mathbf{u}}^0)))\\
&+2\hat{C}_k(\tilde{\mathbf{C}}(\hat{\mathbf{u}}^0,\hat{\mathbf{u}}^0),\tilde{\mathbf{C}}(\hat{\mathbf{u}}^0,\hat{\mathbf{u}}^0)).
\end{align*}Combining these two results gives us the $O(t^2)$ term of the complete memory approximation:
\begin{align}
R_k^2(\hat{\mathbf{u}}) =& Pe^{t\mathcal{L}}[2P\mathcal{L}P\mathcal{L}Q\mathcal{L}-P\mathcal{L}\mathcal{L}Q\mathcal{L}]\hat{u}_k^0\notag\\=&-2\hat{C}_k(\tilde{\mathbf{C}}(\hat{\mathbf{u}},\hat{\mathbf{u}}),\tilde{\mathbf{C}}(\hat{\mathbf{u}},\hat{\mathbf{u}}))\notag\\
&-2\hat{C}_k(\hat{\mathbf{u}},i\epsilon^2 \tilde{\mathbf{C}}^{k3}(\hat{\mathbf{u}},\hat{\mathbf{u}})+2\tilde{\mathbf{C}}(\hat{\mathbf{u}},\tilde{\mathbf{C}}(\hat{\mathbf{u}},\hat{\mathbf{u}}))\notag\\
&\qquad \qquad -2\tilde{\mathbf{C}}(\hat{\mathbf{u}},i\epsilon^2 \hat{\mathbf{u}}^{k3}+\hat{\mathbf{C}}(\hat{\mathbf{u}},\hat{\mathbf{u}}))).
\end{align}One can clearly see how frustrating it would be to compute the $O(t^3)$ and $O(t^4)$ terms by hand. Instead, we use an automated procedure written in Mathematica to derive these expressions. We have reproduced them here:

\begin{align}
R_k^3(\hat{\mathbf{u}}) =& 2\hat{C}_k(\hat{\mathbf{u}},-\epsilon^4 \tilde{\mathbf{C}}^{k6}(\hat{\mathbf{u}},\hat{\mathbf{u}})-2i\epsilon^2\tilde{\mathbf{C}}^{k3}(\hat{\mathbf{u}},2i\epsilon^2 \hat{\mathbf{u}}^{k3}+2\hat{\mathbf{C}}(\hat{\mathbf{u}},\hat{\mathbf{u}})-\tilde{\mathbf{C}}(\hat{\mathbf{u}},\hat{\mathbf{u}}))\notag\\
&\quad \quad \;\;  +2\tilde{\mathbf{C}}(\hat{\mathbf{u}},-\epsilon^4 \hat{\mathbf{u}}^{k6}+i\epsilon^2 \hat{\mathbf{C}}^{k3}(\hat{\mathbf{u}},\hat{\mathbf{u}})+2\hat{\mathbf{C}}(\hat{\mathbf{u}},i\epsilon^2 \hat{\mathbf{u}}^{k3}+\hat{\mathbf{C}}(\hat{\mathbf{u}},\hat{\mathbf{u}})-2\tilde{\mathbf{C}}(\hat{\mathbf{u}},\hat{\mathbf{u}}))\notag\\
&\qquad\qquad\qquad\;+i\epsilon^2 \tilde{\mathbf{C}}^{k3}(\hat{\mathbf{u}},\hat{\mathbf{u}})+2\tilde{\mathbf{C}}(\hat{\mathbf{u}},-2(i\epsilon^2 \hat{\mathbf{u}}^{k3}+\hat{\mathbf{C}}(\hat{\mathbf{u}},\hat{\mathbf{u}}))+\tilde{\mathbf{C}}(\hat{\mathbf{u}},\hat{\mathbf{u}})))\notag\\
&\quad \quad \;\;  +2\tilde{\mathbf{C}}(i\epsilon^2 \hat{\mathbf{u}}^{k3}+\hat{\mathbf{C}}(\hat{\mathbf{u}},\hat{\mathbf{u}}),i\epsilon^2 \hat{\mathbf{u}}^{k3}+\hat{\mathbf{C}}(\hat{\mathbf{u}},\hat{\mathbf{u}})-\tilde{\mathbf{C}}(\hat{\mathbf{u}},\hat{\mathbf{u}}))\notag\\
&\quad \quad \;\;  +2\tilde{\mathbf{C}}(\tilde{\mathbf{C}}(\hat{\mathbf{u}},\hat{\mathbf{u}}),\tilde{\mathbf{C}}(\hat{\mathbf{u}},\hat{\mathbf{u}})))\notag\\
&+6\hat{C}_k(\tilde{\mathbf{C}}(\hat{\mathbf{u}},\hat{\mathbf{u}}),i\epsilon^2  \tilde{\mathbf{C}}^{k3}(\hat{\mathbf{u}},\hat{\mathbf{u}})-2\tilde{\mathbf{C}}(\hat{\mathbf{u}},i\epsilon^2 \hat{\mathbf{u}}^{k3}+\hat{\mathbf{C}}(\hat{\mathbf{u}},\hat{\mathbf{u}})-\tilde{\mathbf{C}}(\hat{\mathbf{u}},\hat{\mathbf{u}}))).
\end{align}
\begin{align}
R_k^4(\hat{\mathbf{u}}) = & 2\hat{C}_k(\hat{\mathbf{u}},i\epsilon^6 \tilde{\mathbf{C}}^{k9}(\hat{\mathbf{u}},\hat{\mathbf{u}})-2\epsilon^4 \tilde{\mathbf{C}}^{k6}(\hat{\mathbf{u}},3i\epsilon^2 \hat{\mathbf{u}}^{k3}+3\hat{\mathbf{C}}(\hat{\mathbf{u}},\hat{\mathbf{u}})-\tilde{\mathbf{C}}(\hat{\mathbf{u}},\hat{\mathbf{u}}))\notag\\
&\qquad\quad-2i\epsilon^2 \tilde{\mathbf{C}}^{k3}(\hat{\mathbf{u}},-3\epsilon^4 \hat{\mathbf{u}}^{k6}+3i\epsilon^2 \hat{\mathbf{C}}^{k3}(\hat{\mathbf{u}},\hat{\mathbf{u}})+2\hat{\mathbf{C}}(\hat{\mathbf{u}},3i\epsilon^2 \hat{\mathbf{u}}^{k3}+3\hat{\mathbf{C}}(\hat{\mathbf{u}},\hat{\mathbf{u}})-5\tilde{\mathbf{C}}(\hat{\mathbf{u}},\hat{\mathbf{u}}))\notag\\
&\qquad\qquad \qquad\qquad\qquad +i\epsilon^2 \tilde{\mathbf{C}}^{k3}(\hat{\mathbf{u}},\hat{\mathbf{u}})-2\tilde{\mathbf{C}}(\hat{\mathbf{u}},3i\epsilon^2 \hat{\mathbf{u}}^{k3}+3\hat{\mathbf{C}}(\hat{\mathbf{u}},\hat{\mathbf{u}})-\tilde{\mathbf{C}}(\hat{\mathbf{u}},\hat{\mathbf{u}})))\notag\\
&\qquad\quad-2\tilde{\mathbf{C}}(\hat{\mathbf{u}},i\epsilon^6 \hat{\mathbf{u}}^{k9}+\epsilon^4 \hat{\mathbf{C}}^{k6}(\hat{\mathbf{u}},\hat{\mathbf{u}})-2i\epsilon^2 \hat{\mathbf{C}}^{k3}(\hat{\mathbf{u}},i\epsilon^2 \hat{\mathbf{u}}^{k3}+\hat{\mathbf{C}}(\hat{\mathbf{u}},\hat{\mathbf{u}})-3\tilde{\mathbf{C}}(\hat{\mathbf{u}},\hat{\mathbf{u}}))\notag\\
&\qquad\qquad\qquad \quad+2\hat{\mathbf{C}}(\hat{\mathbf{u}},\epsilon^4 \hat{\mathbf{u}}^{k6}-i\epsilon^2 \hat{\mathbf{C}}^{k3}(\hat{\mathbf{u}},\hat{\mathbf{u}})-2\hat{\mathbf{C}}(\hat{\mathbf{u}},i\epsilon^2 \hat{\mathbf{u}}^{k3}+\hat{\mathbf{C}}(\hat{\mathbf{u}},\hat{\mathbf{u}})-3\tilde{\mathbf{C}}(\hat{\mathbf{u}},\hat{\mathbf{u}}))\notag\\
&\qquad\qquad\qquad\qquad\qquad\quad -3i\epsilon^2 \tilde{\mathbf{C}}^{k3}(\hat{\mathbf{u}},\hat{\mathbf{u}})+2\tilde{\mathbf{C}}(\hat{\mathbf{u}},5i\epsilon^2 \hat{\mathbf{u}}^{k3}+5\hat{\mathbf{C}}(\hat{\mathbf{u}},\hat{\mathbf{u}})-3\tilde{\mathbf{C}}(\hat{\mathbf{u}},\hat{\mathbf{u}})))\notag\\
&\qquad\qquad\qquad \quad-2\hat{\mathbf{C}}(i\epsilon^2 \hat{\mathbf{u}}^{k3} +\hat{\mathbf{C}}(\hat{\mathbf{u}},\hat{\mathbf{u}}),i\epsilon^2 \hat{\mathbf{u}}^{k3}+\hat{\mathbf{C}}(\hat{\mathbf{u}},\hat{\mathbf{u}})-2\tilde{\mathbf{C}}(\hat{\mathbf{u}},\hat{\mathbf{u}}))\notag\\
&\qquad\qquad\qquad \quad-6\hat{\mathbf{C}}(\tilde{\mathbf{C}}(\hat{\mathbf{u}},\hat{\mathbf{u}}),\tilde{\mathbf{C}}(\hat{\mathbf{u}},\hat{\mathbf{u}}))\notag\\
&\qquad\qquad\qquad \quad-\epsilon^4 \tilde{\mathbf{C}}^{k6}(\hat{\mathbf{u}},\hat{\mathbf{u}})+2i\epsilon^2 \tilde{\mathbf{C}}^{k3}(\hat{\mathbf{u}},-3i\epsilon^2 \hat{\mathbf{u}}^{k3}-3\hat{\mathbf{C}}(\hat{\mathbf{u}},\hat{\mathbf{u}})+\tilde{\mathbf{C}}(\hat{\mathbf{u}},\hat{\mathbf{u}}))\notag\\
&\qquad\qquad\qquad \quad+2\tilde{\mathbf{C}}(\hat{\mathbf{u}},-3\epsilon^4 \hat{\mathbf{u}}^{k6}+3i\epsilon^2 \hat{\mathbf{C}}^{k3}(\hat{\mathbf{u}},\hat{\mathbf{u}})+2\hat{\mathbf{C}}(\hat{\mathbf{u}},3i\epsilon^2 \hat{\mathbf{u}}^{k3}+3\hat{\mathbf{C}}(\hat{\mathbf{u}},\hat{\mathbf{u}})-5\tilde{\mathbf{C}}(\hat{\mathbf{u}},\hat{\mathbf{u}}))\notag\\
&\qquad\qquad\qquad\qquad\qquad\quad +i\epsilon^2 \tilde{\mathbf{C}}^{k3}(\hat{\mathbf{u}},\hat{\mathbf{u}})+2\tilde{\mathbf{C}}(\hat{\mathbf{u}},-3i\epsilon^2 \hat{\mathbf{u}}^{k3}-3\hat{\mathbf{C}}(\hat{\mathbf{u}},\hat{\mathbf{u}})+\tilde{\mathbf{C}}(\hat{\mathbf{u}},\hat{\mathbf{u}})))\notag\\
&\qquad\qquad\qquad \quad+2\tilde{\mathbf{C}}(i\epsilon^2 \hat{\mathbf{u}}^{k3} +\hat{\mathbf{C}}(\hat{\mathbf{u}},\hat{\mathbf{u}}),3i\epsilon^2 \hat{\mathbf{u}}^{k3}+3\hat{\mathbf{C}}(\hat{\mathbf{u}},\hat{\mathbf{u}})-2\tilde{\mathbf{C}}(\hat{\mathbf{u}},\hat{\mathbf{u}}))\notag\\
&\qquad\qquad\qquad \quad+2\tilde{\mathbf{C}}(\tilde{\mathbf{C}}(\hat{\mathbf{u}},\hat{\mathbf{u}}),\tilde{\mathbf{C}}(\hat{\mathbf{u}},\hat{\mathbf{u}})))\notag\\
&\qquad\quad+2i\epsilon^2 \tilde{\mathbf{C}}^{k3}(i\epsilon^2 \hat{\mathbf{u}}^{k3}+\hat{\mathbf{C}}(\hat{\mathbf{u}},\hat{\mathbf{u}}),-3\epsilon^2\hat{\mathbf{u}}^{k3}-3\hat{\mathbf{C}}(\hat{\mathbf{u}},\hat{\mathbf{u}})+2\tilde{\mathbf{C}}(\hat{\mathbf{u}},\hat{\mathbf{u}}))\notag\\
&\qquad\quad-2\tilde{\mathbf{C}}(i\epsilon^2 \hat{\mathbf{u}}^{k3}+\hat{\mathbf{C}}(\hat{\mathbf{u}},\hat{\mathbf{u}}),3\epsilon^4 \hat{\mathbf{u}}^{k6}-3i\epsilon^2 \hat{\mathbf{C}}^{k3}(\hat{\mathbf{u}},\hat{\mathbf{u}})+2\hat{\mathbf{C}}(\hat{\mathbf{u}},-3i\epsilon^2 \hat{\mathbf{u}}^{k3}-3\hat{\mathbf{C}}(\hat{\mathbf{u}},\hat{\mathbf{u}})+5\tilde{\mathbf{C}}(\hat{\mathbf{u}},\hat{\mathbf{u}}))\notag\\
&\qquad\qquad\qquad \qquad\qquad\qquad\quad\quad-i\epsilon^2 \tilde{\mathbf{C}}^{k3}(\hat{\mathbf{u}},\hat{\mathbf{u}})+2\tilde{\mathbf{C}}(\hat{\mathbf{u}},3i\epsilon^2 \hat{\mathbf{u}}^{k3}+3\hat{\mathbf{C}}(\hat{\mathbf{u}},\hat{\mathbf{u}})-\tilde{\mathbf{C}}(\hat{\mathbf{u}},\hat{\mathbf{u}})))\notag\\
&\qquad\quad+2\tilde{\mathbf{C}}(\tilde{\mathbf{C}}(\hat{\mathbf{u}},\hat{\mathbf{u}}),\epsilon^4 \hat{\mathbf{u}}^{k6}-i\epsilon^2 \hat{\mathbf{C}}^{k3}(\hat{\mathbf{u}},\hat{\mathbf{u}})-2\mathbf{C}(\hat{\mathbf{u}},i\epsilon^2 \hat{\mathbf{u}}^{k3}+\hat{\mathbf{C}}(\hat{\mathbf{u}},\hat{\mathbf{u}})-3\tilde{\mathbf{C}}(\hat{\mathbf{u}},\hat{\mathbf{u}}))\notag\\
&\qquad\qquad\qquad\qquad\qquad\;\;-3i\epsilon^2 \tilde{\mathbf{C}}^{k3}(\hat{\mathbf{u}},\hat{\mathbf{u}})+2\tilde{\mathbf{C}}(\hat{\mathbf{u}},5i\epsilon^2 \hat{\mathbf{u}}^{k3}+5\hat{\mathbf{C}}(\hat{\mathbf{u}},\hat{\mathbf{u}})-3\tilde{\mathbf{C}}(\hat{\mathbf{u}},\hat{\mathbf{u}})))\notag\\
&\qquad\quad-2i\epsilon^2 \tilde{\mathbf{C}}^{k3}(\tilde{\mathbf{C}}(\hat{\mathbf{u}},\hat{\mathbf{u}}),\tilde{\mathbf{C}}(\hat{\mathbf{u}},\hat{\mathbf{u}})))\notag\\
&-8\hat{C}_k(\tilde{\mathbf{C}}(\hat{\mathbf{u}},\hat{\mathbf{u}}),-\epsilon^4 \tilde{\mathbf{C}}^{k6}(\hat{\mathbf{u}},\hat{\mathbf{u}})+2i\epsilon^2 \tilde{\mathbf{C}}^{k3}(\hat{\mathbf{u}},-2i\epsilon^2 \hat{\mathbf{u}}^{k3}-2\hat{\mathbf{C}}(\hat{\mathbf{u}},\hat{\mathbf{u}})+\tilde{\mathbf{C}}(\hat{\mathbf{u}},\hat{\mathbf{u}}))\notag\\
&\qquad\qquad \qquad\;\;\;\;\;\; +2\tilde{\mathbf{C}}(\hat{\mathbf{u}},-\epsilon^4 \hat{\mathbf{u}}^{k6}+i\epsilon^2 \hat{\mathbf{C}}^{k3}(\hat{\mathbf{u}},\hat{\mathbf{u}})+2\hat{\mathbf{C}}(\hat{\mathbf{u}},i\epsilon^2 \hat{\mathbf{u}}^{k3}+\hat{\mathbf{C}}(\hat{\mathbf{u}},\hat{\mathbf{u}})-2\tilde{\mathbf{C}}(\hat{\mathbf{u}},\hat{\mathbf{u}}))\notag\\
&\qquad\qquad \qquad\qquad \qquad\qquad +i\epsilon^2 \tilde{\mathbf{C}}^{k3}(\hat{\mathbf{u}},\hat{\mathbf{u}})+2\tilde{\mathbf{C}}(\hat{\mathbf{u}},-2i\epsilon^2 \hat{\mathbf{u}}^{k3}-2\hat{\mathbf{C}}(\hat{\mathbf{u}},\hat{\mathbf{u}})+\tilde{\mathbf{C}}(\hat{\mathbf{u}},\hat{\mathbf{u}})))\notag\\
&\qquad\qquad \qquad\;\;\;\;\;\; +2\tilde{\mathbf{C}}(i\epsilon^2 \hat{\mathbf{u}}^{k3}+\hat{\mathbf{C}}(\hat{\mathbf{u}},\hat{\mathbf{u}}),i\epsilon^2 \hat{\mathbf{u}}^{k3}+\hat{\mathbf{C}}(\hat{\mathbf{u}},\hat{\mathbf{u}})-\tilde{\mathbf{C}}(\hat{\mathbf{u}},\hat{\mathbf{u}}))\notag\\
&\qquad\qquad \qquad\;\;\;\;\;\; +2\tilde{\mathbf{C}}(\tilde{\mathbf{C}}(\hat{\mathbf{u}},\hat{\mathbf{u}}),\tilde{\mathbf{C}}(\hat{\mathbf{u}},\hat{\mathbf{u}})))\notag\\
&+48\hat{C}_k(\tilde{\mathbf{C}}(\hat{\mathbf{u}},i\epsilon^2 \hat{\mathbf{u}}^{k3}+\hat{\mathbf{C}}(\hat{\mathbf{u}},\hat{\mathbf{u}})),i\epsilon^2 \tilde{\mathbf{C}}^{k3}(\hat{\mathbf{u}},\hat{\mathbf{u}})+2\tilde{\mathbf{C}}(\hat{\mathbf{u}},\tilde{\mathbf{C}}(\hat{\mathbf{u}},\hat{\mathbf{u}})))\notag\\
&-6\hat{C}_k(i\epsilon^2 \tilde{\mathbf{C}}^{k3}(\hat{\mathbf{u}},\hat{\mathbf{u}})+2\tilde{\mathbf{C}}(\hat{\mathbf{u}},i\epsilon^2 \hat{\mathbf{u}}^{k3}+\hat{\mathbf{C}}(\hat{\mathbf{u}},\hat{\mathbf{u}})+\tilde{\mathbf{C}}(\hat{\mathbf{u}},\hat{\mathbf{u}})),\notag\\
&\qquad\qquad\qquad\qquad \qquad i\epsilon^2 \tilde{\mathbf{C}}^{k3}(\hat{\mathbf{u}},\hat{\mathbf{u}})+2\tilde{\mathbf{C}}(\hat{\mathbf{u}},i\epsilon^2 \hat{\mathbf{u}}^{k3}+\hat{\mathbf{C}}(\hat{\mathbf{u}},\hat{\mathbf{u}})+\tilde{\mathbf{C}}(\hat{\mathbf{u}},\hat{\mathbf{u}})))
\end{align}
\bibliography{MultiScale}
\bibliographystyle{plain}

\end{document}